\documentclass[12pt]{amsart}
\usepackage{amsmath, amsthm, amssymb, amsfonts}
\usepackage[english]{babel}
\usepackage{geometry}
\usepackage{xcolor}
\usepackage{enumitem}
\usepackage{bbm}
\usepackage{indentfirst}
\usepackage{amsaddr}

\theoremstyle{definition}
\newtheorem{lemma}{Lemma}
\newtheorem{theorem}[lemma]{Theorem}
\newtheorem{definition}[lemma]{Definition}
\newtheorem{proposition}[lemma]{Proposition}
\newtheorem{corollary}[lemma]{Corollary}

\newtheorem{question}[lemma]{Question}
\DeclareMathOperator{\andd}{and}

\DeclareMathOperator{\cf}{cf}
\DeclareMathOperator{\Fin}{Fin\times Fin}
\newcommand{\conc}{^\smallfrown}

\DeclareMathOperator{\dom}{dom}
\DeclareMathOperator{\Pseudo}{Pseudo}
\newcommand{\rest}{\upharpoonright} 

\def\cA{{\mathcal{A}}} \def\cB{{\mathcal{B}}}   \def\cE{{\mathcal{E}}} \def\cF{{\mathcal{F}}}   \def\cI{{\mathcal{I}}} \def\cJ{{\mathcal{J}}}         \def\cS{{\mathcal{S}}}  \def\cV{{\mathcal{V}}} \def\cU{{\mathcal{U}}}    
\title{Fin-intersecting MAD families}

\author[C. Corral and V. O. Rodrigues]{Cesar Corral and Vinicius de O. Rodrigues}

\address[C. Corral, V. O. Rodrigues]{York University}
\email{cicorral@yorku.ca and vor@yorku.ca}
\makeatletter
\renewcommand{\email}[2][]{%
	\ifx\emails\@empty\relax\else{\g@addto@macro\emails{,\space}}\fi%
	\@ifnotempty{#1}{\g@addto@macro\emails{\textrm{(#1)}\space}}%
	\g@addto@macro\emails{#2}%
}
\makeatother
\begin{document}
	\maketitle
	
	\begin{abstract}
		We introduce a new class of almost disjoint families which we call fin-intersecting almost disjoint families. They are related to almost disjoint families whose Vietoris Hyperspace of their Isbell-Mrówka spaces are pseudocompact. We show that under $\mathfrak p=\mathfrak c$ fin-intersecting MAD families exist generically and they also exist if $\mathfrak{a<s}$, but that there are also non fin-intersecting MAD families in ZFC. We also show that under CH, there exists fin intersecting MAD families which remain like that after adding an arbitrary quantity of Cohen reals and Random reals. These results give more models in which pseudocompact MAD families exist.
	\end{abstract}
	\section{Introduction}

	\subsection{Some history}
	Topologies on spaces of subsets of a given topological space have been studied since the beginning of the last century. In this paper we aim our attention to the Vietoris Hyperspace which is defined as follows:
	
	\begin{definition}
		Let $X$ be a $T_1$ Topological space. By $\exp(X)$ we denote the set of all nonempty closed subsets of $X$.
		
		Given a subset $U$ of $X$, we let $U^+=\{F\in \exp(X):F\subseteq U\}$ and  $U^-=\{F\in \exp(X):F\cap U\neq \emptyset\}$.
		
		The \textit{Vietoris topology} on $\exp(X)$ is the topology generated by $\{U^+, U^-: U\subseteq X\text{ is open}\}$.
		
		The \textit{Vietoris hyperspace} of $X$ is the set $\exp(X)$ endowed with the Vietoris topology.\qed
	\end{definition}
	
	Some of the first steps towards topologizing collections of subsets of a given topological
	space were taken by Hausdorff \cite{Hausdorff1927}, who defined a metric on $\exp(X)$
	in the case where $X$ is a bounded metric space. This metric is usually called Hausdorff's metric. The Vietoris Topology was introduced by Leopold Vietoris \cite{vietoris1922bereiche} and coincides with the topology generated by the Hausdorff's metric in case $X$ is a compact metric space.
	
	Although $X$ does not need to be $T_1$ for this definition to make sense, we restrict ourselves to $T_1$ spaces to avoid some pathologies.
	
	The topological properties of $X$ often depends on the topological properties of $\exp(X)$ and vice-versa. For instance, $X$ is normal iff $\exp(X)$ is regular iff $\exp(X)$ is Tychonoff, and $X$ is regular iff $\exp(X)$ is Hausdorff. A nice reference for the basics on the Vietoris topology is \cite{michael1951topologies}.
	
	Leopold Vietoris has proved that $X$ is compact iff its Vietoris hyperspace is compact \cite{vietoris1922bereiche}. This result motivates the question of whether there are relations between $X$ and $\exp(X)$ concerning generalizations of compactness. These questions were investigated by J. Ginsburg, who proved some results on this direction regarding pseudocompactness, countably compactness, $p$-compactness and $p$-pseudocompactness (where $p$ is some fixed free ultrafilter) \cite{ginsburg1975some}. In particular, he proved that if every power of $X$ is countably compact, so does $\exp(X)$, but he could not prove the same for pseudocompactness.
	
	Later, M. Hru\v{s}\'{a}k, I. Martínes-Ruiz and F. Hernandez-Hernandez proved that there exists a subspace $X$ of $\beta \omega$ containing $\omega$ such that $X^\omega$ is pseudocompact but $\exp(X)$ is not pseudocompact (recall that $X^\omega$ is pseudocompact iff every power of $X$ is pseudocompact) \cite{hrusak2007pseudocompactness}. This example was later improved by V. Rodrigues, A. Tomita and Y. Ortiz-Castillo by making all powers of $X$ less than the cardinal characteristic $\mathfrak h$ countably compact \cite{ORTIZCASTILLO20189}. 
	
	These examples illustrate that it is difficult to infer the pseudocompactness of $\exp(X)$ by knowing topological properties of the powers of $X$, so it makes sense to study this problem restricted to particular classes of spaces. A natural class of spaces to explore this problem are the Isbell-Mrówka spaces, which are topological spaces associated to almost disjoint families. An almost disjoint family on a countable infinite set $N$ is an infinite set $\mathcal A$ of infinite subsets of $N$ which are pairwise almost disjoint, that is, for every $a, b\in \mathcal A$, $a\cap b$ is finite. A maximal almost disjoint family (on $N$), also called MAD family, is an almost disjoint family (on $N$) which is not contained in any other almost disjoint family. MAD families exist by Zorn's lemma and are uncountable. The least size of a MAD family is denoted by $\mathfrak a$. By $\mathcal I(\mathcal A)$ we denote the free ideal generated by $\mathcal A$, that is, $\mathcal I(\mathcal A)=\{X\subseteq \omega:\exists \mathcal B\in [\mathcal A]^{<\omega}\,X\subseteq^*\bigcup \mathcal B\}$.

	Given an almost disjoint family $\mathcal A$ one defines its \textit{Isbell-Mrówka space}, also called its \textit{psi space}, as follows:
	
	\begin{definition}
		Let $\mathcal A$ be an almost disjoint family on a countable infinite set $N$ such that $N\cap [N]^\omega=\emptyset$ (such as $\omega$ or $2^{\omega}$). $\Psi(\mathcal A)$ is the set $\omega\cup\mathcal A$ topologized by the finest topology which makes $\omega$ open and discrete and each $a\subseteq \Psi(\mathcal A)$ is a sequence converging to the point $a\in \Psi(\mathcal A)$. A basis for this topology is $\{\{n\}:n\in N\}\cup\{\{a\}\cup(a\setminus F):a\in \mathcal A\andd F\in [N]^{<\omega}\}$.
	\end{definition}
	
	Of course, the set $N$ is combiratorically and topologically unimportant if no extra structure is being discussed. It is straightforward to see that $\Psi(\mathcal A)$ is a Hausdorff, locally compact (therefore Tychonoff) zero-dimensional separable first-countable not countably compact topological space.
	
	An excellent survey on this class of spaces is \cite{Hrusak2014}. The topological properties of $\Psi(\mathcal A)$ often depend on the combinatorical properties of $\mathcal A$.  The following is what makes psi spaces interesting for examining the pseudocompactness of hyperspaces of Vietoris:
	
	\begin{proposition}
		Let $\mathcal A$ be an almost disjoint family. Then $\mathcal A$ is a MAD family iff $\Psi(\mathcal A)$ is pseudocompact iff $\Psi(\mathcal A)^\omega$ is pseudocompact.\qed
	\end{proposition}
	
	As a shorthand, we say that an almost disjoint family $\mathcal A$ is pseudocompact iff $\exp(\Psi(\mathcal A))$ is pseudocompact.
	
	It is known that if $\exp(X)$ is pseudocompact, so is $X$ \cite{ginsburg1975some}. Thus, when restricting to psi spaces, the questions regarding the relations between the pseudocompactness of $X$ and of $\exp(X)$ boil down to the following:
	
	\begin{question}\label{question1}Is there a pseudocompact MAD family?
	\end{question}
	\begin{question}\label{question2}Is every MAD family pseudocompact?
	\end{question}
	
	Question \ref{question2} has been partially answered in \cite{hrusak2007pseudocompactness} and was later solved in \cite{NewPaperPseudocompactness}, where we proved that every MAD family is pseudocompact iff the Baire number of $\omega^*$ is strictly greater than $\mathfrak c$. Question \ref{question1} appears in \cite{Hrusak2014}, \cite{hruvsak2018pseudocompact} and \cite{hrusak2007pseudocompactness} and is still open.
	
	We do not have many examples of non pseudocompact MAD families (or even of pseudocompact ones in case $n(\omega^*)\leq\mathfrak c$). In particular, the non pseudocompact MAD family constructed in \cite{hrusak2007pseudocompactness} using $\mathfrak h<\mathfrak c$ has cardinality $\mathfrak c>\omega_1$, and we constructed one of size $\omega_2<\mathfrak c$ in \cite{NewPaperPseudocompactness}. It is unknown whether it is consistent to exist a non pseudocompact MAD family of size $\omega_1$.
	
	To further study Question \ref{question2} we introduce, in this paper, a new class of almost disjoint families which we call fin-intersecting. Fin-intersecting MAD families are pseudocompact.
	
	\subsection{Structure of the paper and summary of results}

	In \cite{guzman2017generic}, O. Guzman and M. Hru\v{s}\'{a}k have defined the notion of generic existence of an almost disjoint family with a property $\phi$. Given a property $\phi$ about almost disjoint families, they have defined that \textit{almost disjoint families with the property $\phi$} exist generically iff every (infinite) almost disjoint family of cardinality $<\mathfrak c$ may be extended to an almost disjoint family with the property $\phi$. So, for instance, MAD families exist generically iff $\mathfrak a=\mathfrak c$. Note that this concept has nothing to do with forcing or generic filters and it  is much more linked to the idea of extending a MAD family by recursion to achieve a certain property.
	
	\begin{itemize}
		\item In Section 2 we introduce fin-intersecting almost disjoint families and show that fin-intersecting MAD families are pseudocompact.
		\item In Section 3 we show that every almost disjoint family of size $<\mathfrak s$ is fin-intersecting. In particular, if $\mathfrak a<\mathfrak c$ there is a fin-intersecting MAD family. We show that fin-intersecting MAD families exist generically under $\mathfrak{ap=s=c}$. These results give, in particular, examples of more models containing pseudocompact MAD families, partially answering Question \ref{question2}. Moreover, it is easy to observe that fin-intersecting almost disjoint families do not need to be MAD, but it is non-trivial to see if they must be MAD somewhere. Thus, we show that assuming $\mathfrak{b=s=c}$ fin-intersecting nowhere MAD families exist generically.
		\item In Section 4, we show that there exists a non fin-intersecting MAD family in ZFC. The proof is divided in two cases: $\mathfrak h<\mathfrak c$ and $\mathfrak s\leq \mathfrak a$. This shows that (consistently) not every pseudocompact MAD family is fin-intersecting. We also ask if $\mathfrak s$ is the smallest cardinal for which there is a non fin-intersecting almost-disjoint family and show that this is the case if $\mathfrak s\leq \mathfrak {ie}$.
		
		\item In Section 5 we show that by assuming CH there exist fin-intersecting MAD families which remain fin-intersecting MAD after adding any number of Cohen/Random reals using finite supports. Note that in these models $\mathfrak s=\omega_1$ so these are non-trivial results. This also shows that pseudocompact MAD families exist in the Random model, what was not known.

	\end{itemize}
	
	For the sake of completeness we define the cardinals mentioned above:
	
	\begin{itemize}
		\item $\mathfrak{ap}$ is the least size of an almost disjoint family $\mathcal A$ for which there exists $\mathcal B\subseteq \mathcal A$ such that there is no $c\in [\omega]^\omega$ such that for every $a\in \mathcal A$ $c\cap a$ is infinite iff $a\in \mathcal B$.
		\item $\mathfrak{b}$ is the least size of a collection $\mathcal B\subseteq \omega^\omega$ such that for every $g\in \omega^\omega$ there exists $f\in \mathcal B$ such that $f\not \leq^* g$, where $f\leq^*g$ means that $\{n \in \omega: g(n)<f(n)\}$ is finite.
		\item $\mathfrak  h$ is the least size of a collection of open dense subsets of $[\omega]^\omega$ whose intersection is empty (equivalently: not dense). We say that $\mathcal D\subseteq [\omega]^\omega$ is dense iff for every $a\in [\omega]^\omega$ there exists $b\in \mathcal D$ such that $b\subseteq^* a$, and it is open iff for every $b \in \mathcal D$ and  $a\in [\omega]^\omega$, if $a\subseteq^* b$ then $a\in \mathcal D$.
		\item To define $\mathfrak p$, we say that a collection $\mathcal P\subseteq [\omega]^\omega$ is centered iff the intersection of every finite nonempty subcollection of $\mathcal P$ is infinite. A pseudointersection for $\mathcal P$ is an infinite subset $P\subseteq \omega$ such that $P\subseteq^*A$ for every $A\in \mathcal P$. It is clear that if $\mathcal P$ admits a pseudointersection then it is centered. The cardinal $\mathfrak p$ is the least cardinality of a centered collection $\mathcal P\subseteq [\omega]^\omega$  with no pseudointersection.
		\item $\mathfrak s$ is the least cardinality of a collection $\mathcal S\subseteq [\omega]^\omega$ such that for every $a\subseteq \omega$ there exists $b\in \mathcal S$ such that $|a\cap b|=\omega$ and $|a\setminus b|=\omega$. Such an $\mathcal S$ is called a splitting family.
		\item To define $\mathfrak{ie}$, let $\triangle=\{(n, m)\in \omega\times\omega: m\leq n\}$. Then $\mathfrak{ie}$ is the least size of a collection of partial infinite functions $\mathcal F$ from $\omega$ into $\omega$ contained in $\triangle$ such that no single total function from $\omega$ into $\omega$ contained in $\triangle$ is almost disjoint from every element of $\mathcal F$.
		
	\end{itemize}
	
	For more on $\mathfrak{ap}$ see \cite{Brendle}. For more on $\mathfrak{a, b, h, s}$ see \cite{Blass2010}, and $\mathfrak{ie}$ was defined in \cite{NewPaperPseudocompactness}.
	
	In the references above one may find proofs for $\omega_1\leq \mathfrak p\leq \mathfrak h, \mathfrak {ap}\leq \mathfrak b, \mathfrak s\leq \mathfrak c$, for $\mathfrak b\leq \mathfrak a\leq \mathfrak c$. A routine diagonalization argument shows that $\omega_1\leq\mathfrak{ie}\leq \mathfrak c$.
	
    Moreover, it is well known that $\mathfrak p\leq \cf(\mathfrak s)$. For the sake of completeness we quickly sketch a proof: if $(\kappa_\alpha: \alpha<\theta)$ is a increasing sequence of cardinals less than $\mathfrak s$ indexed by some cardinal $\theta<\mathfrak p$ and $\lambda=\sup_{\alpha<\mathfrak p}\kappa_\alpha$, then $\lambda<\mathfrak s$ as well. To see this, if $\mathcal S=\{S_\xi:\xi<\lambda\}\subseteq [\omega]^\omega$, one easily recursively constructs a $\subseteq^*$-decreasing tower $b_\alpha\in [\omega]^\omega$ for $\alpha<\theta$ such that for every $\xi<\kappa_\alpha$, $b_{\alpha+1}\subseteq S_\xi$ or $b_{\alpha+1}\cap S_\xi=^*\emptyset$. Now a pseudointersection $b$ of this tower shows that $\mathcal S$ is not a splitting family. In particular, this shows that $\cf(\mathfrak s)>\omega$.

	\section{Fin-intersecting almost disjoint families}
	
	We shall adopt the following definition which is helpful while exploring pseudocompactness in Isbell-Mrówka spaces.
	\begin{definition}
		A fin sequence is a function $C:\omega\rightarrow [\omega]^{<\omega}\setminus \{\emptyset\}$ such that for all $n, m \in \omega$, if $n\neq m$ then $C(n)\cap C(m)=\emptyset$.
	\end{definition}
	
	In \cite[Proposition 2.1.]{Rodrigues2019small}, V. Rodrigues and A. Tomita proved the following theorem:
	
	\begin{proposition}Let $\mathcal A$ be an almost disjoint family. Then $\mathcal A$ is pseudocompact iff every fin sequence has an accumulation point in $\exp(\Psi(\mathcal A))$.
	\end{proposition}
	
	The nice thing about the previous proposition is that we may restrict ourselves with working with pairwise disjoint families. We are ready to introduce fin-intersecting almost disjoint families, which are closely related with fin sequences and pseudocompactness. Recall that a collection of sets $\mathcal P$ is centered if every finite nonempty subcollection of $\mathcal P$ has infinite intersection.
	
	\begin{definition}We say that an almost disjoint family $\mathcal A$ is fin-intersecting iff for every fin sequence $C$ there exists an infinite set $I\subseteq \omega$ such that $\{\{n \in I: a \cap C_n \neq \emptyset\}: a \in \mathcal A\}\setminus [I]^{<\omega}$ is centered.
		
	\end{definition}
	
	Thus, in order of an almost disjoint family $\mathcal A$ be fin-intersecting for every fin sequence $C$ we need to have an infinite set $I$ such that after throwing away every finite set of $\{\{n\in I: a\cap C_n\neq \emptyset\}: a\in \mathcal A\}$, the remaining set is centered.

	The next lemma helps us to show that certain sequences automatically satisfy what is required in the definition of fin-intersecting almost disjoint families. 
	
	\begin{lemma}Let $\mathcal A$ be an almost disjoint family and $C$ be a fin sequence. Suppose there exists an infinite $I\subseteq \omega$ such that $\bigcup_{n \in I} C_n \in \mathcal I(\mathcal A)$. Then there exists $J\in[I]^\omega$ such that $\{\{n\in J: C(n)\cap a\neq \emptyset\}:a\in\mathcal A\}\setminus [J]^{<\omega}$ is centered.\label{lemma:idealImpliesFinIntersecting}
	\end{lemma}
	
	\begin{proof}
		Let $k$ be the least natural number for which there exists $\mathcal B\subseteq \mathcal A$ of size $k$ and $J\in[I]^\omega$ such that $\bigcup_{n\in J}C(n)\subseteq^*\bigcup \mathcal B$. Then:
		
		\begin{enumerate}[label=\arabic*)]
			\item For every $a\in \mathcal B$, $\{n\in J:C(n)\cap a=\emptyset\}$ is finite, for if it was infinite, $J'=\{n\in J:C(n)\cap a=\emptyset\}$ and $\mathcal B\setminus\{a\}$ would violate the minimality of $k$.
			\item For every $a\in \mathcal A\setminus \mathcal B$, $\{n\in J: C(n)\cap a\neq\emptyset\}$ is finite, for if it was infinite we would have $a\cap\bigcup_{n\in J}C(n)\subseteq^*\bigcup \mathcal B$, which implies that for some $b\in \mathcal B$, $a\cap b$ is infinite, a contradiction.
		\end{enumerate}
		
		Thus $\{\{n\in J: C(n)\cap a\neq \emptyset\}:a\in\mathcal A\}\setminus [J]^{<\omega}=\{\{n\in J: C(n)\cap a\neq \emptyset\}:a\in\mathcal B\}$, and the intersection of the latter set is a cofinite subset of $J$ by 1).
	\end{proof}
	
	Now we show that fin-intersecting MAD families are pseudocompact, which is our main motivation for studying them.
	\begin{proposition}
		Every fin-intersecting MAD family is pseudocompact.
	\end{proposition}
	\begin{proof}
		Let $\mathcal A$ be a fin-intersecting MAD family and $C$ a fin-sequence. Fix $I$ as in the definition of fin-intersecting almost disjoint families. Let $\mathcal B=\{a \in \mathcal B: |\{n \in I: C(n) \cap a\neq \emptyset \}|=\omega\}$. Notice that $\mathcal B$ is a nonempty (and closed) set since $\mathcal A$ is MAD and the $C_n$'s are pairwise disjoint. We claim that $\mathcal B$ is an accumulation point for $C$.

		Let $V^+\cap U_0^-\cap \dots\cap U_m^-$ be a basic open neighborhood of $\mathcal B$, where $V, U_0, \dots, U_m$ are open subsets of $\Psi(\mathcal A)$. By intersecting each $U_i$ with $V$, we may assume that for each $i\leq m$, $U_i\subseteq V$.	
		
		For each $i\leq m$ there exists $a \in \mathcal B$ and $N_i \in \omega$ such that $\{a_i\}\cup (a_i\setminus N_i) \subseteq U_i$. Let $J=\bigcap_{i \leq m}\{n \in I: a_i \cap C_n \neq \emptyset\}\setminus \{n\in I: C_n\cap \bigcup_{i\leq m}N_i=\emptyset\}$, which is infinite. Notice that if $n \in J$, then for each $i\leq m$, $C_n\cap a_i\setminus N_i\neq \emptyset$ so $C_n\cap U_i\neq \emptyset$. We claim that for all but finitely many $n \in J$, $C_n\subseteq V$, which completes the proof. 
		
		Suppose that this is not true. Then $J'=\{n \in J: C_n\setminus  V\neq \emptyset\}$ is infinite. Since the $C_n$'s are pairwise disjoint, $X=\bigcup_{n \in J'}C_n\setminus V$ is infinite. Since $\mathcal A$ is MAD, there exists $a \in \mathcal A$ such that $|a\cap X|=\omega$.
		
		\textbf{Case 1:} $a \in \mathcal A\setminus\mathcal B$. In this case, $\{n \in I: C_n\cap a\neq \emptyset\}$ is finite, so $a\cap \bigcup_{n \in I} C_n\supseteq a\cap X$ is finite, a contradiction.
		
		\textbf{Case 2:} $a \in \mathcal B$. In this case, $a\subseteq^*V$, so $|a\cap X|<\omega$, a contradiction.
	\end{proof}
	
	Thus, many questions concerning fin intersecting almost disjoint families automatically appear. Is every pseudocompact MAD family fin intersecting? Are all fin-intersecting almost disjoint families MAD (somewhere)? Do fin-intersecting MAD families exist?
	
	We will partly or completely answer all these questions. Also, fin-intersecting MAD families may be easier to work with than pseudocompact MAD families, so we were able to use them to produce some new examples of pseudocompact MAD families.
	\section{On the existence of fin-intersecting almost disjoint families}
	
	In this section we prove that under certain conditions fin-intersecting MAD families exist. Some of these results imply the existence of pseudocompact MAD families in contexts where their existence was not known. We also show that under certain conditions fin-intersecting nowhere MAD families exist, which proves that they may not need to be somewhere MAD.

	The cardinal invariant $\mathfrak{s}$ is closely tied to fin-intersecting almost disjoint families. Recall that for $S,X\in[\omega]^\omega$ we say that $S$ splits $X$ if $|X\cap S|=\omega=|X\setminus S|$. A splitting family $\cS\subseteq[\omega]^\omega$ is a family such that for every $X\in[\omega]^\omega$ there exists $S\in \cS$ such that $S$ splits $X$. Then $\mathfrak{s}$ is the minimum cardinality of a splitting family.
	
	\begin{theorem}\label{<simpliesfinintersecting}
		Every almost disjoint family of size less than $\mathfrak{s}$ is fin-intersecting.
	\end{theorem}
	
	\begin{proof}
		Let $\cA=\{a_\alpha:\alpha<\kappa\}$ be an almost disjoint family with $\kappa<\mathfrak{s}$. Given a fin sequence $C:\omega\to[\omega]^{<\omega}\setminus\emptyset$ define $S_\alpha=\{n\in\omega:a_\alpha\cap C(n)\neq\emptyset\}$. Thus since $\cS=\{S_\alpha:\alpha<\kappa\}$ is not a splitting family, there exists $I\in[\omega]^\omega$ such that either $I\subseteq^* S_\alpha$ or $I\cap S_\alpha=^*\emptyset$. Now it is clear that for every $\cF\in [\cA]^{<\omega}$, if $\{n\in\omega:a\cap C(n)\neq\emptyset\}$ is infinite for every $a\in\mathcal F$, then $I\subseteq^*\{n\in I:\forall\ a\in\cF\, (a\cap C(n)\neq\emptyset)\}$.
	\end{proof}

	\begin{corollary}\label{a<s}
		If $\mathfrak{a<s}$ there exists a fin-intersecting MAD family. In particular, there are pseudocompact MAD families if $\mathfrak{a<s}$.
	\end{corollary}
	
	So there are pseudocompact MAD families if $\mathfrak a<\mathfrak s$. This was previously unknown.
	
	Recall that under $\mathfrak p=\mathfrak c$ every MAD family is pseudocompact (\cite{hrusak2007pseudocompactness}, \cite{NewPaperPseudocompactness}). $\mathfrak p=\mathfrak c$ also implies the existence of many fin-intersecting MAD families, as illustrated by the result below. Recall that $\mathfrak p\leq\mathfrak{ap}\leq \mathfrak{a}$ \cite{Brendle} and that $\mathfrak p\leq  \mathfrak s$ \cite{Blass2010}.
	
	We will need the following result:
	
	\begin{theorem}[{\cite[Theorem 2.1.]{mateweak}}]Let $(a_i: i\in I)$ be a family of infinite subsets of $\omega$ such that $|I|<\mathfrak c$. Then there exists an indexed almost disjoint family $(b_i: i\in I)$ such that $b_i\subseteq a_i$ for every $i \in I$.
	\end{theorem}
	
	\begin{theorem}[$\mathfrak {ap}=\mathfrak c=\mathfrak s$]\label{genexistencefinintersect} Fin-intersecting MAD families exist generically.
	\end{theorem}

	\begin{proof}
		Let $\kappa<\mathfrak c$ be an infinite cardinal and $\mathcal A=\{a_\alpha:\alpha<\kappa\}$ be an almost disjoint family enumerated injectively. We must extend $\mathcal A$ to a fin-intersecting MAD family.
		
		Enumerate $[\omega]^\omega=\{X_\alpha:\kappa\leq\alpha<\mathfrak c\}$ and all fin sequences as $\{C_\alpha:\kappa\leq\alpha<\mathfrak c\}$. We recursively define $I_\alpha$ and $a_\alpha$ for $\alpha\in[\kappa, \mathfrak c)$ satisfying that for every $\alpha\in [\kappa, \mathfrak c)$:
		\begin{enumerate}
			\item $a_\alpha\in [\omega]^\omega$ and for every $\beta<\alpha$, $a_\beta\cap a_\alpha$ is finite.
			
			\item $I_\alpha\in [\omega]^\omega$ and for every $\beta<\alpha$, either $\{n\in I_\alpha: C_\alpha(n)\cap a_\beta=\emptyset\}$ or $\{n\in I_\alpha: C_\alpha(n)\cap a_\beta\neq \emptyset\}$ are finite.
			\item For every $\gamma \in [\kappa, \alpha]$, and for every $F \in [\alpha]^{<\omega}$,  \\
			
			if for all $J \in [\alpha]^{<\omega}$ $|\{n \in I_\alpha: \forall \xi \in F\, (C_\gamma(n)\cap a_\xi\neq \emptyset) \andd  C_\gamma(n)\setminus \bigcup_{\xi \in J}{a_\xi}\neq \emptyset\}|=\omega$,\\
			
			then $|\{n \in I_\alpha: \forall \xi \in F\cup \{\alpha\}\,(C_\gamma(n)\cap a_\xi\neq \emptyset) \}|=\omega$.
			
			\item If for every $\beta<\alpha$ $X_\alpha\cap a_\beta$ is finite, then $|a_\alpha\cap X_\alpha|=\omega$.
		\end{enumerate}

		At step $\alpha$ we first define $I_\alpha$ as follows: let $\mathcal S_\alpha=\{\{n\in \omega: C_\alpha(n)\cap a_\beta=\emptyset\}: \beta<\alpha\}$. Since $|\mathcal S_\alpha|<\mathfrak c=\mathfrak s$, $\mathcal S_\alpha$ is not a splitting family, so there exists $I_\alpha\in [\omega]^\omega$ such that for every $\beta<\alpha$, $\{n\in I_\alpha: C_\alpha(n)\cap a_\beta=\emptyset\}$ is finite or $\{n\in I_\alpha: C_\alpha(n)\cap a_\beta\neq\emptyset\}$ is finite, as intended.
		
		Now we define $a_\alpha^0$ satisfying 4. and 1, $a_\alpha^1$ satisfying 3. and 1., so obviously $a_\alpha=a_\alpha^0\cup a_\alpha^1$ will satisfy 1., 3. and 4., which finishes the construction.
		
		Defining $a_\alpha^0$: if the hypothesis of 4. does not hold, just let $a_\alpha^0$ satisfy 1., which is possible since $\mathfrak a=\mathfrak c$. If the hypothesis hold just let $a_\alpha^0=X_\alpha$.
		
		Defining $a_\alpha^1$: if the hypothesis of 3. does not hold, just let $a_\alpha^1$ satisfy 1., which is possible since $\mathfrak a=\mathfrak c$. Now suppose it holds. For each $\gamma\in [\kappa,\alpha]$ and $F, J\in [\alpha]^{<\omega}$ let:
		
		$$K_\gamma(F, J)=\left\{n\in I_\gamma: \forall \xi \in F \,(C_\gamma(n)\cap a_\xi\neq\emptyset)\andd C_\gamma(n)\setminus \bigcup_{\xi\in J}a_\xi\neq\emptyset\right\}, \andd$$
		
		$$K_\gamma(F)=\left\{n\in I_\gamma: \forall \xi \in F \,(C_\gamma(n)\cap a_\xi\neq\emptyset)\right\}.$$
		
		For each $\gamma \in [\kappa, \mathfrak c)$, let $\mathcal F_\gamma=\{F\in [\alpha]^{<\omega}:\forall J\in [\alpha]^{<\omega}\, |K_\gamma(F, J)|=\omega\}$. Now for each $F \in \mathcal F_\gamma$, since for every $J \in [\alpha]^{<\omega}$ we have $\bigcup_{n\in K_\gamma(F)}C(n)\in \mathcal I^+(\{a_\xi:\xi<\alpha\})$. By using $\mathfrak a=\mathfrak c$ there exists $a_\gamma(F)\subseteq \bigcup_{n\in K_\gamma(F)}C(n)$ such that $a_\gamma(F)\cap a_\xi$ is finite for every $\xi<\alpha$. Since the family $(a_\gamma(F): \gamma \in [\kappa, \alpha], F\in \mathcal F_\gamma)$ has size less than $\mathfrak c$, we may shrink the elements to an indexed almost disjoint family $\cU=(b_\gamma(F): \gamma \in [\kappa, \alpha], F\in \mathcal F_\gamma)$.
		
		Let $\mathcal V=\{a_\xi:\xi<\alpha\}$. Notice that $|\mathcal U|,|\mathcal V|<\mathfrak c=\mathfrak{ap}$ and $\cU\cup\cV$ is almost disjoint. By $\mathfrak {ap}=\mathfrak c$, there exists $a^1_\alpha\in [\omega]^\omega$ such that $a^1_\alpha \cap x$ is infinite for every $x\in \mathcal V$ and $a^1_\alpha\cap a_\xi$ is finite for every $\xi<\alpha$. To verify 3., given $\gamma \in [\kappa, \alpha]$, $F\in [\alpha]^{<\omega}$ such that for every $J\in [\alpha]^{<\omega}$, $|\{n \in I_\alpha: \forall \xi \in F\, (C_\gamma(n)\cap a_\xi\neq \emptyset) \andd  C_\gamma(n)\setminus \bigcup_{\xi \in J}{a_\xi}\neq \emptyset\}|=\omega$, we have that $F\in \mathcal F_\gamma$, so $a^1_\alpha\cap b_\gamma(F)$ is infinite, which implies that $a_\alpha^1\cap \bigcup_{K_\gamma(F)}C_\gamma(n)$ is also infinite and this implies that $\{n \in I_\alpha: \forall \xi \in F\cup \{\alpha\}\,(C_\gamma(n)\cap a_\xi\neq \emptyset) \}=\{n\in K_\gamma(F): C_\gamma(n)\cap a_\alpha^1\neq\emptyset\}$ is infinite.\\
		
		Now we verify that $\{a_\alpha:\alpha<\mathfrak c\}$ is a fin-intersecting MAD family.
		
		1. and 4. easily imply that it is a MAD family.
		
		To see that it is fin-intersecting, let $C$ be a fin sequence. In case there exists an infinite $I\subseteq \omega$ such that $\bigcup_{n\in I}C_n\in \mathcal I(\{a_\alpha:\alpha<\mathfrak c\})$ we are done by Lemma \ref{lemma:idealImpliesFinIntersecting}, so assume this does not happen.
		
		There exists $\gamma$ such that $C_\gamma=C$. We show that $I=I_\gamma$ works as in the definition of a fin-intersecting almost disjoint family. So let $F\in [\mathfrak c]^{<\omega}$ be nonempty and such that for every $\xi \in F$, $\{n\in I_\gamma: C_\gamma(n)\cap \xi\}$ is infinite. We must see that $\bigcap_{\xi\in F}\{n\in I:C_\gamma(n)\cap a_\xi\neq\emptyset\}$ is infinite. Let $F_0=F\cap\gamma$, write $F\setminus F_0=\{\xi_0, \dots, \xi_{m-1}\}$ in increasing order and define $F_{i+1}=F_i\cup \{\xi_i\}$ for $1\leq i\leq m$. We show by finite induction that for each $i\leq m$, $\bigcap_{\xi\in F_i}\{n\in I:C_\gamma(n)\cap a_\xi\neq\emptyset\}$ is infinite.
		
		By 2., $\bigcap_{\xi\in F_0}\{n\in I:C_\gamma(n)\cap a_\xi\neq\emptyset\}$ is cofinite in $I_\gamma$ (in case $F_0=\emptyset$ we consider this empty intersection to be $I_\gamma$). Now suppose that $i<m$ and that we know that $L=\bigcap_{\xi\in F_i}\{n\in I_\gamma:C_\gamma(n)\cap a_\xi\neq\emptyset\}$ is infinite. Then for every finite $J\in [\xi_i]^{<\omega}$ we know that $\bigcup_{n\in L}C_\gamma(n)\setminus \bigcup_{\xi\in J}a_\xi$ is infinite, which implies that $\{n\in L: C_\gamma(n)\setminus \bigcup_{\xi\in J}a_\xi\neq\emptyset \}$ is infinite, but this set is $\{n \in I_\gamma: \forall \xi \in F_i\, (C_\gamma(n)\cap a_\xi\neq \emptyset) \andd  C_\gamma(n)\setminus \bigcup_{\xi \in J}{a_\xi}\neq \emptyset\}$. By 2. with $\alpha=\xi_i$, this implies that $\{n \in I_\gamma: \forall \xi \in F_i\cup \{\xi_i\}\,(C_\gamma(n)\cap a_\xi\neq \emptyset) \}$ is infinite, but this set is $\bigcap_{\xi\in F_{i+1}}\{n\in I_\gamma:C_\gamma(n)\cap a_\xi\neq\emptyset\}$. This completes the proof.
	\end{proof}
	
	The previous theorem shows that under certain cardinal hypothesis many fin-intersecting MAD families exist. Since these MAD families are pseudocompact, this result potentially yields more results about the existence of pseudocompact MAD families as well. However, under $\mathfrak p=\mathfrak c$ (or, more weakly, under $n(\omega^*)>\mathfrak c$, which is false under $\mathfrak h<\mathfrak c$) all MAD families are pseudocompact \cite{NewPaperPseudocompactness}. Thus, the following question interests us:
	
	\begin{question}
		Is $\mathfrak{ap}=\mathfrak s=\mathfrak c\geq n(\omega^*)$ consistent?
	\end{question}
	
	It is natural to ask whether there is any relation between madness and being fin-intersecting. From Theorem \ref{<simpliesfinintersecting} we know that not every fin-intersecting almost disjoint family is maximal and that consistently, there is an uncountable fin-intersecting almost disjoint family that is not MAD (e.g., any non maximal almost disjoint family of size $\omega_1$ in a model of $\mathfrak{s=c>\omega_1}$). Moreover, since the fin-intersecting property is preserved under subfamilies, if $\cA$ is a fin-intersecting MAD family and $a_0\in\cA$, the subfamily $\cA\setminus\{a_0\}$ is still fin-intersecting and fails to be MAD. However the previous example is MAD when restricted to $\omega\setminus a_0$. Thus, it is natural to ask if fin-intersecting almost disjoint families of size $\mathfrak c$ are always somewhere MAD.
	
	To address that question we show that fin-intersecting almost disjoint families that fail to be MAD everywhere (consistently) exist generically. Despite the almost disjoint family in Theorem \ref{genexistencefinintersect} being MAD, the construction of every $a_\alpha$ is split into two parts, $a_\alpha=a_\alpha^0\cup a_\alpha^1$, where $a_\alpha^1$ witnesses the fin-intersecting property and $a_\alpha^0$ is a witness for madness. We can modify the latter construction of $a_\alpha^0$ and define $a_\alpha=a_\alpha^1$ in order to get a nowhere MAD almost disjoint family.

	Recall that the the ideal $\Fin$ is the ideal on $\omega\times\omega$ consisting of the subsets $A$ of $\omega\times \omega$ such that there exists a function $f:\omega\rightarrow \omega$ and a cofinite set $B\subseteq \omega$ such that for every $(n, m)\in A$, if $n\in B$ then $m\leq f(n)$ (this is the Fubini product of the ideal of the finite subsets of $\omega$ by itself).

	\begin{theorem}[$\mathfrak{b}=\mathfrak{c}=\mathfrak{s}$]
		Nowhere MAD fin-intersecting almost disjoint families exist generically.
	\end{theorem}
	
	\begin{proof}
		Let $\mathcal A$ be an almost disjoint family of size $\kappa<\mathfrak c$. We may enumerate it as $\mathcal A=\{a_\alpha: \alpha<\kappa\}$. It is easy to see that fin-intersecting nowhere MAD families are preserved by finite modifications in the members of the family, so we may assume that $\{A_n:n\in\omega\}$ forms a partition and moreover, without loss of generality, $\mathcal A$ is an almost disjoint family on $\omega\times \omega$ and for every $n\in \omega$, $a_n=\{n\}\times \omega$.
		
		Let $[\omega\times \omega]^\omega=\{X_\alpha: \kappa\leq \alpha<\mathfrak c\}$ and let $(C_\alpha: \kappa\leq \alpha<\mathfrak c)$ be the family of all fin-sequences $C:\omega\rightarrow [\omega\times\omega]^{<\omega}\setminus\{\emptyset\}$.
		
		For $\alpha \in [\kappa, \mathfrak c)$ we recursively define sets $I_\alpha \in [\omega]^\omega$ and infinite sets $a_\alpha, b_\alpha \in [\omega\times\omega]^\omega$ such that:
		
		\begin{enumerate}[label=\alph*)]
			\item For every $\xi<\alpha$, $|a_\xi\cap a_\alpha|<\omega$ and $|a_\xi\cap b_\alpha|<\omega$.
			\item For every $\xi\in[\kappa, \alpha]$, $|b_\xi\cap a_\alpha|<\omega$.
			
			\item For every $\xi\in[\kappa,\alpha]$, $I_\alpha\subseteq^* \{n \in \omega: C_\alpha(n)\cap a_\xi\neq \emptyset\}$ or $I_\alpha\subseteq^* \{n \in \omega: C_\alpha(n)\cap a_\xi=\emptyset\}$.
			\item For every $\gamma \in [\kappa, \alpha]$, if $\bigcup_{n\in I_\gamma}C_\gamma(n) \in \Fin$ then $\{n \in I_\gamma: C_\gamma(n)\cap a_\alpha\neq \emptyset\}$ is finite.
			\item For every $\gamma \in [\kappa, \alpha]$, if $\bigcup_{n\in I_\gamma}C_\gamma(n) \notin \Fin$ then for every $F\in [\alpha]^{<\omega}$,
			$$\big|\{n \in I_\gamma: \forall \xi \in F\,( C_\gamma(n)\cap a_\xi\neq \emptyset)\}\big|=\omega\implies$$
			$$|\{n \in I_\gamma: \forall \xi \in (F\cup\{\alpha\})( C_\gamma(n)\cap a_\xi\neq \emptyset)\}|=\omega$$
			\item If $X_\alpha \in \mathcal I^+(\{a_\xi: \xi<\alpha\})$ then $b_\alpha \subseteq X_\alpha$.
			
		\end{enumerate}
		
		We show that this construction is possible. Suppose we have constructed everything up to step $\alpha<\mathfrak c$.
		
		We fix $I_\alpha$ satisfying c) which is possible since $\alpha<\mathfrak c=\mathfrak s$.
		
		Now we define $b_\alpha$. If the hypothesis of f) fails, just let $b_\alpha$ be almost disjoint with every member of $\{a_\xi: \xi<\alpha\}$ which is possible since $\mathfrak a=\mathfrak c$. If it holds, $\{a_\xi\cap X_\alpha:\xi<\alpha\}$ is not MAD since $\mathfrak a=\mathfrak c$, then there exists $b_\alpha$ as in f) and a).
		
		It remains to define $a_\alpha$. First, notice that for every $\xi\in [\omega, \alpha)$ and for every $\eta \in  [\kappa, \alpha]$, $a_\xi$ and $b_\beta$ have finite intersection with each column $a_n$, so we fix:
		
		\begin{itemize}
			\item for every $\xi\in [\kappa, \alpha)$ let $f_\xi\in \omega^\omega$ be such that $a_\xi\subseteq\{(n, m)\in \omega\times\omega: m< f_\xi(n)\}$.
			\item for every $\eta\in [\kappa, \alpha)$ let $g_\eta\in \omega^\omega$ be such that $b_\eta\subseteq\{(n, m)\in \omega\times\omega: m<g_\eta(n)\}$.
			
			\item $\Gamma=\{\gamma\in [\kappa, \alpha]:\bigcup_{n\in I_\gamma}C_\gamma(n)\in \Fin \}$
			
			\item for every $\gamma\in\Gamma$, let $h_\gamma\in \omega^\omega$ be such that there exists a finite $J\subseteq \omega$ such that $\bigcup_{n\in I_\gamma}C_\gamma(n)\subseteq\{(n, m)\in \omega\times\omega: m< h_\gamma(n)\}\bigcup (J\times \omega)$.
			
		\end{itemize}
		
		Since $\mathfrak b=\mathfrak c$ there exists $f \in \omega^\omega$ such that for every $\xi\in [\omega, \alpha)$ and for every $\gamma\in \Gamma$, we have $f_\xi, g_\xi, h_\gamma\leq^*f$.
		
		It is clear that any infinite subset of $\omega\times \omega$ bounded by a function $F$ and that is above $f$ will satisfy a), b), d). It remains to select $F$ that makes e) hold.
		
		Let $\gamma\in[\kappa,\alpha]$ and $E\in[\alpha]^{<\omega}$. Assume that $X(\gamma)=\bigcup_{n\in I_\gamma} C_\gamma(n)\notin\Fin$ and let $Z(\gamma,\xi)=\{n\in I_\gamma:a_\xi\cap C_\gamma(n)\neq\emptyset\}$ for every $\xi<\alpha$.
		From c), d) and e), it is clear that if $Z(\gamma,\eta)$ is infinite for every $\eta\in G\in[\alpha]^{<\omega}$, then $\bigcap_{\eta\in G}Z(\gamma,\eta)$ is infinite too.
		
		Assume then that $Z(\gamma,\xi)$ is infinite for every $\xi\in E$. Since $X(\gamma)\notin\Fin$, there are infinitely many $n\in\omega$ such that $Z(\gamma,n)$ is infinite. Let $\{k_i:i\in\omega\}$ be the increasing enumeration of this set. Define simultaneously $F(\gamma,E)\in\omega^\omega$ and an increasing sequence $\{n_i:i\in\omega\}$ as follows:
		
		If $k\notin\{k_i:i\in\omega\}$ let $F(\gamma,E)(k)=f(k)$ and for every $k_i$, find $n_i$ and define $F(\gamma,E)(k_i)$ such that $\left(k_i\times[f(k_i),F(\gamma,E)(k_i)]\cap C_\gamma(n_i)\neq\right)\emptyset$ and $n_i\in \bigcap_{\eta\in E}Z(\gamma,\eta)$. This choice is possible since $\bigcap_{\eta\in E\cup\{k_i\}}Z(\gamma,\eta)$ is infinite for every $k_i$. 
		
		Let $F\in\omega^\omega$ be such that $F\geq^* F(\gamma,E)$ for every $\gamma\in[\kappa,\alpha]$ and $E\in[\alpha]^{<\omega}$. Define $a_\alpha=\bigcup_{n\in\omega}\{n\}\times[f(n),F(n)]$. Since $\bigcup_{n\in\omega}\{n\}\times[f(n),F(\gamma,E)(n)]$ $\subseteq^* a_\alpha$ for every pair $(\gamma,E)$, item e) holds. 
		
		Finally define $\cA=\{a_\alpha:\alpha<\mathfrak{c}\}$ and notice that for every $X\in\cI^+(\cA)$, if $X=X_\alpha$, then $X_\alpha\in\cI^+(\{a_\beta:\beta<\alpha\})$ implies that $b_\alpha\subseteq X$, that is, $\cA$ is nowhere MAD. It is also easy to show inductively like in Theorem \ref{genexistencefinintersect} that $\cA$ is fin-intersecting.
	\end{proof}

	
		
		We do not know if the cardinal hypothesis for the two previous theorems on fin-intersecting (nowhere) MAD families can be weakened or if they even exist in \textsf{ZFC}.
		
		\begin{question}
			Is there a fin-intersecting MAD family in $\textsf{ZFC}$?
		\end{question}
		
		\begin{question}
			Is the existence of a fin-intersecting MAD family consistent with $\mathfrak{s<a}$?
		\end{question}
		
		\begin{question}
			Is there a nowhere MAD fin-intersecting family of size $\mathfrak{c}$ in \textsf{ZFC}?
		\end{question}
		
		\section{Non fin-intersecting almost disjoint families}
		
		One could ask if under $\mathfrak p=\mathfrak c$ every MAD family is fin-intersecting, just like it happens with pseudocompact MAD families. The answer is negative: we will see that non fin-intersecting MAD families exist in ZFC. In particular, this shows, as we shall see, that pseudocompact MAD family which are not fin-intersecting consistently exist. Thus, fin-intersecting madness is not a characterization of pseudocompactness.
		
		Given a fin-sequence $C:\omega\to[\omega]^{<\omega}$, a set $X\in[\omega]^\omega$ is a \emph{selector} for $C$ if $|C(n)\cap X|\leq1$ for every $n\in\omega$.
		We say that a fin sequence $C$ is \emph{unbounded} if
		$$|\{n\in\omega:|C(n)|<k\}|<\omega$$
		for every $k\in\omega$.
		
		Recall that a MAD family is completely separable iff for every $X \in \mathcal I^+(\mathcal A)$ there exists $a \in \mathcal A$ such that $a\subseteq X$. Completely separable MAD families were defined by S. Hechler \cite{hechler1971classifying} to study problems related to $\beta \omega$. It is unknown if they exist in ZFC, but it is known that they exist in many models of set theory \cite{Hrusak2014}.
		
		\begin{lemma}
			Let $C$ be an unbounded fin-sequence. Then every completely separable MAD family consisting of selectors for $C$ is not fin-intersecting.
		\end{lemma}
		
		\begin{proof}
			Let $\cA$ be a mad family consisting of selectors for $C$ and let $I\in[\omega]^\omega$. Let $I=I_0\sqcup I_1$ be a partition. Since $\cA$ consists of selectors, $X=\bigcup_{n\in I_0}C(n)$ and $Y=\bigcup_{n\in I_1}C(n)$ are disjoint sets in $\cI^+(\cA)$. 
			Thus there exist $a_0,a_1\in\cA$ such that $a_0\subseteq X$ and $a_1\subseteq Y$. Hence $\{n\in I:a_i\cap C(n)\neq\emptyset\}=I_i$ and in consequence
			$$\{\{n\in I:a\cap C(n)\neq\emptyset\}:a\in\cA\}\setminus[I]^{<\omega}$$
			is not centered. Since $I$ was chosen arbitrarily, $\cA$ is not fin-intersecting.
		\end{proof}
		
		We now show that families satisfying the previous hypothesis consistently exist. The following theorem is an easy modification of the construction of a completely separable MAD family under $\mathfrak{s}<\mathfrak{a}$ given by Shelah \cite{shelah2011mad} and later improved by Mildenberger, Raghavan and Steprans, weakening the hypothesis to $\mathfrak{s}\leq\mathfrak{a}$ \cite{mildenberger2014splitting}. We include its proof for the sake of completeness.
		
		\begin{theorem}[$\mathfrak{s}\leq\mathfrak{a}$]\label{specialcompsepmad}
			Let $C$ be a fin-sequence (not necessarily unbounded), then there is a completely separable MAD family $\cA$ consisting of selectors for $C$.
		\end{theorem}
		
		\begin{proof}
			Let $C$ be a fin sequence. For $X\subseteq\omega$ let $X^0=X$ and $X^1=\omega\setminus X$.
			Let $\{S_\alpha:\alpha<\mathfrak{s}\}$ be a family such that for every almost disjoint family $\cB$ and every $X\in\cI^+(\cB)$ there exists $\alpha<\mathfrak{s}$ such that $X\cap S_\alpha,\ X\setminus S_\alpha\in\cI^+(\cB)$ (the existence of this family was proven in \cite[Lemma 7]{raghavan2012weakly} for $\mathfrak s_{\omega,\omega}$ in the place of $\mathfrak s$, and in \cite[Theorem 6]{mildenberger2014splitting} it was shown that $\mathfrak s=\mathfrak s_{\omega, \omega}$).\\
			
			Enumerate $[\omega]^\omega=\{X_\alpha:\alpha<\mathfrak{c}\}$.
			We recursively construct $\cA=\{a_\alpha:\alpha<\mathfrak{c}\}$ and $\{\tau_\alpha:\alpha<\mathfrak{c}\}\subseteq2^{<\mathfrak{s}}$ such that for every $\xi<\alpha<\mathfrak{c}$
			\begin{enumerate}
				\item $|a_\alpha\cap a_\xi|<\omega$,
				\item $\tau_\alpha\nsubseteq\tau_\xi$,
				\item $a_\alpha\subseteq^* S_\eta^{\tau_\alpha(\eta)}$ for every $\eta\in\dom(\tau_\alpha)$,
				\item $|a_\alpha\cap C(n)|\leq1$ for every $n\in\omega$ and
				\item if $X_\alpha\in\cI^+(\{a_\beta:\beta<\alpha\})$ then $a_\alpha\subseteq X_\alpha$.
				\item $\bigcup_{n\in \omega}C_n\setminus \bigcup_{i<k}a_k$ is infinite for every $k\in \omega$ (we need this to guarantee that the construction does not end before $\omega$).
			\end{enumerate}
			
			At step $\alpha$, assume without loss of generality (by (6)) that $X_\alpha\in\cI^+(\{a_\beta:\beta<\alpha\})$. Define $X_\emptyset=X$ and for every $t\in2^{<\omega}$ define a function $\sigma_t$, $X_{t\conc0}\in \mathcal I^+(\{a_\beta: \beta<\alpha\}$ and $X_{t\conc1}\in \mathcal I^+(\{a_\beta: \beta<\alpha\}$ as follows:
			\begin{itemize}
				\item $\dom(\sigma_t)=\min\{\beta<\mathfrak{s}:X_t\cap S_\beta,X_t\setminus S_\beta\in\cI^+(\{a_\beta:\beta<\alpha\})\}$,
				\item $\sigma_t(\eta)$ is the unique $i\in2$ such that $X_t\cap S_\eta^{\sigma_t(\eta)}\in\cI^+(\{a_\beta:\beta<\alpha\})$,
				\item if $\beta=\dom(\sigma_t)$ then $X_{t\conc0}=X_t\cap S_\beta^0\in\cI^+(\{a_\beta:\beta<\alpha\})$ and $X_{t\conc1}=X_t\cap S_\beta^1\in\cI^+(\{a_\beta:\beta<\alpha\})$.
			\end{itemize}
			Once we do this, we end up with families $\{X_t\in2^{<\omega}\}\subseteq\cI^+(a_\beta:\beta<\alpha)$ and $\sigma_t$ such that $s\subseteq t$ implies $X_t\subseteq X_s$, $\sigma_s\subset\sigma_t$ and if $s,t$ are incomparable, then so are $\sigma_s$ and $\sigma_t$ and $X_t\cap X_s=\emptyset$. Recall that given an almost disjoint family $\cB$ and decreasing sequence $(Y_n:n\in\omega)$ in $\cI^+(\cB)$, there exists a pseudointersection $Y\in\cI^+(\cB)$ (see for instance \cite[Proposition 2.3]{Hrusak2014}). Thus for every $f\in2^\omega$ we can find $X_f\in\cI^+(\{a_\beta:\beta<\alpha\})$ pseudointersection for $\{X_{f\rest n}:n\in\omega\}$ and define also $\sigma_f=\bigcup_{n\in\omega}\sigma_{f\rest n}$ (recall that $\mathfrak{s}$ has uncountable cofinality).
			
			It is clear that $\{\sigma_f:f\in2^\omega\}$ consists on incomparable elements, so we can find $f\in2^\omega$ such that no $\tau_\xi$ extends $\sigma_f$. Define $\tau_\alpha=\sigma_f$ and let $\delta=\dom(\tau_\alpha)$. Let $\cB=\{a_\beta:\beta<\alpha\land \tau_\beta\subseteq\tau_\alpha\}$. Also, for every $\eta<\delta$ let $\cF_\eta\in[\{a_\beta:\beta<\alpha\}]^{<\omega}$ be such that $X_f\setminus S_\eta^{\tau_\alpha(\eta)}\subseteq^*\bigcup\cF_\eta$. Then $\cE=\Big|\cB\cup\bigcup\big\{\cF_\eta:\eta<\delta\big\}\Big|<\mathfrak{s\leq a}$. Let $a_\alpha\subseteq X_f$ almost disjoint with every element in $\cE$. By shrinking $a_\alpha$ if necessary, we can assume that 4 and 6 hold, and clearly 2 and 5 hold. It is also easy to see that 1 follows from 3 and the fact that $a_\alpha$ is almost disjoint with any element in $\cB$. Now we show that 3 holds.
			
			Let $\eta<\dom(\tau_\alpha)$, since $a_\alpha\subseteq X_f$ and $X_f\setminus S_\eta^{\tau_\alpha(\eta)}\subseteq^*\bigcup \cF_\eta$, we have that $a_\alpha\setminus S_\eta^{\tau_\alpha(\eta)}\subseteq^*\bigcup\cF_\eta$. On the other hand, $|a_\alpha\cap \bigcup\cF_\eta|<\omega$ and therefore (3) holds.
			
			Now we see that (1) holds. Fix $\xi<\alpha$. If $\tau_\xi\subseteq \tau_\alpha$ this follows from the choice of $\mathcal B$, $\mathcal E$ and $ a_\alpha$. If $\tau_\xi\not\subseteq \tau_\alpha$, since $\tau_\alpha \not \subseteq \tau_\xi$ by construction, there exists $\eta \in \delta \cap \dom(\tau_\alpha)$ 
			such that $\tau_\alpha(\eta)=1-\tau_\xi(\eta)$. By (3) we have $a_\alpha\subseteq^* S_\eta^{\tau_\alpha(\eta)}$ and $a_\xi\subseteq^* S_\eta^{\tau_\alpha(\eta)}$, so $a_\alpha\cap a_\xi\subseteq^* S_\eta^{\tau_\alpha(\eta)}\cap S_\eta^{\tau_\xi(\eta)}=\emptyset$. Thus, (1) holds.
			
			It is now straightforward to check that this is the desired family.
		\end{proof}
		
		Combining this result with the following, and since $\mathfrak{h\leq s,a}$, we get that non fin-intersecting MAD families exists in \textsf{ZFC}.
		
		\begin{theorem}{\cite{hrusak2007pseudocompactness}}\label{hernandezhrusakmartinez}
			($\mathfrak{h<c}$) There exists a non pseudocompact MAD family.    
		\end{theorem}
		
		\begin{corollary}
			There exists a non fin-intersecting MAD family.
		\end{corollary}
		
		\begin{proof}
			If $\mathfrak{h<c}$, there is a non pseudocompact MAD family, and since fin-intersecting MAD families are pseudocompact, that family necessarily fails to be fin-intersecting. On the other hand if $\mathfrak{h=c}$ then $\mathfrak{s=a=c}$ and the conclusion follows from Theorem \ref{specialcompsepmad}.
		\end{proof}
		
		We do not known if the example in \cite{hrusak2007pseudocompactness} can be made completely separable or if at least there is a completely separable non fin-intersecting MAD family under the same assumption.
		
		\begin{question}
			Is there a completely separable MAD family that is not fin-intersecting in \textsf{ZFC}? Equivalently under $\mathfrak{h<c}$? 
		\end{question}
		
		Also, since the proof of the existence of a non fin-intersecting MAD family goes by cases, we do not known of a single example of a MAD family that is not fin-intersecting in \textsf{ZFC}.
		
		\begin{question}
			Is there an explicit example of a single MAD family $\cA$ that is not fin-intersecting in \textsf{ZFC}?
		\end{question}
		
		Under $\mathfrak{p}=\mathfrak{c}$, maximality of almost disjoint families is equivalent to being pseudocompact. Due to the previous result, since $\mathfrak{p}=\mathfrak{c}$ implies $\mathfrak{s}=\mathfrak{a}=\mathfrak{c}$, we get that there is a completely separable MAD family that is pseudocompact but fails to be fin-intersecting. 
		
		\begin{proposition}[$\mathfrak p=\mathfrak c$] There exists a completely separable pseudocompact MAD family that is not fin-intersecting.
		\end{proposition}

		We do not know if the result above can be improved to the Baire number of $\omega^*$ being greater than $\mathfrak c$.\\
		
		We now look into a different question: Theorem \ref{<simpliesfinintersecting} says that every almost disjoint family of size $<\mathfrak s$ is fin-intersecting. We conjecture that this result is the best possible we can achieve in the sense that we conjecture that $\mathfrak s$ is the first cardinal $\kappa$ with a non fin-intersecting almost disjoint family of size $\kappa$. We will see that this is true under $\mathfrak s\leq \mathfrak {ie}$.

			
			
			
		
		\begin{proposition}
			If $\mathfrak{s\leq ie}$, there exists a non fin-intersecting almost disjoint family of size $\mathfrak{s}$.
		\end{proposition}
		
		\begin{proof}
			Let $\{S_\alpha:\alpha<\mathfrak{s}\}$ be a splitting family. We   recursively construct $\cA=\{a_\alpha^0,a_\alpha^1:\alpha\in\mathfrak{s}\}$ in $\triangle$ such that for every $\alpha<\mathfrak{s}$:
			
			\begin{enumerate}
				\item $\{a_\beta^i:\beta<\alpha\land i\in2\}$ is almost disjoint,
				\item $a_\beta^i$ is a partial function with domain $S_\alpha^i$ where $S_\alpha^0=S$ and $S_\alpha^1=\omega\setminus S_\alpha$.
			\end{enumerate}
			
			At step $\alpha<\mathfrak{s}$, since $\{a_\beta^i:\beta<\alpha\land i\in2\}$ has size less that $\mathfrak{s\leq ie}$, there exists a total function $a_\alpha$ in $\triangle$ almost disjoint with every $a_\beta^i$. Let $a_\alpha^i=a\rest S_\alpha^i$.
			
			To prove that it is not fin-intersecting let $C:\omega\to [\triangle]^{<\omega}$ given by $C(n)=\{(n,m):m\leq n\}$. Thus 
			$$\{n\in\omega:a_\alpha^i\cap C(n)\neq\emptyset\}=S_\alpha^i$$ 
			for every $\alpha<\mathfrak s$ and $i<2$. Let $I\in[\omega]^\omega$. There exists $\alpha<\mathfrak{s}$ such that $|I\cap S_\alpha|=\omega=|I\setminus S_\alpha|$.
			In this case, $\{n\in I: a_\alpha^i\cap C(n)\neq\emptyset\}=I\cap S_\alpha^i$, both are infinite but 
			$$\{n\in I:a_\alpha^0\cap C(n)\neq\emptyset\neq a_\alpha^1\cap C(n)\}=(I\cap S_\alpha^0)\cap(I\cap S_\alpha^1)=\emptyset,$$
			which shows that the family in the definition of fin-intersecting is not centered and thus $\cA$ cannot be fin-intersecting.
		\end{proof}
		
		We do not know if $\mathfrak{s}$ characterize the minimal size of a non fin-intersecting family.
		
		\begin{question}
			Is it true in \textsf{ZFC} that there is a non fin-intersecting almost disjoint family of size $\mathfrak{s}$?
		\end{question}
		
		Also, since the MAD family of Theorem \ref{hernandezhrusakmartinez} and completely separable MAD families have size $\mathfrak{c}$, our examples of non fin-intersecting MAD families are all of size $\mathfrak{c}$. Of course, since every MAD family of size less than $\mathfrak{s}$ is fin-intersecting, the best we can ask for is the following:
		
		\begin{question}
			Does it follow from $\mathfrak{c}>\mathfrak{a\geq s}$ that there are MAD families of size $\mathfrak{a}$ that fail to be fin-intersecting?
		\end{question}
		
		\section{A fin-indestructible MAD family in the Cohen Model}
		
		In this section we show that assuming CH, every countable almost disjoint family can be extended to a fin-intersecting MAD family which remains a fin-intersecting MAD family after adding an arbitrary quantity of Cohen reals (with finite supports). This construction has many similarities to the construction in \cite{Rodrigues2019small} and is a modification of K. Kunen's construction of a Cohen-indestructible MAD family \cite{kunenOld}. The construction in \cite{Rodrigues2019small} inspired the definition of fin-intersecting almost disjoint families.
		
		The following notation comes in handy:
		
		\begin{definition}
			Let $A=(b_n: n \in \omega)$ be a centered countable family of elements of $[\omega]^\omega$. We define $\Pseudo(A)=\{\min(\bigcap_{k\leq n}b_k\setminus n): n\in\omega\}$.
		\end{definition}
		Notice that $\Pseudo(A)$ is a pseudointersection of $\{b_n: n \in \omega\}$ and that $\Pseudo(A)$ is absolute for transitive models of ZFC.
		
		\begin{definition}Let $\beta$ be an infinite countable ordinal, $f:\beta\rightarrow \omega$ be a bijection, $\mathcal A=\{a_\alpha: \alpha<\beta\}$ be an injectively enumerated countable almost disjoint family and $C$ be a fin sequence. We define inductively:
			\begin{itemize}
				\item $I_0(\mathcal A, C, f)=\omega$,
				\item $I_{n+1}=\{i \in I_n: C(i)\cap a_{f^{-1}(n)}\neq \emptyset\}$ if this set is infinite,
				\item $I_{n+1}=\{i \in I_n: C(i)\cap a_{f^{-1}(n)}= \emptyset\}$ if the set above is finite.	
			\end{itemize}
			
			We let $I(\mathcal A, C, f)=\Pseudo((I_n(\mathcal A, C, f))_{n \in \omega})$.
		\end{definition}
		
		Notice that the concepts above are absolute for transitive models of ZFC. The important feature of $I=I(\mathcal A, C, f)=\Pseudo((I_n(\mathcal A, C, f))_{n \in \omega})$ is that it is infinite and that given $a \in \mathcal A$, either $\{i \in I: C(i)\cap a = \emptyset\}$ is finite or $\{i \in I: C(i)\cap a\neq \emptyset\}$ is finite. Thus, $\{\{n \in I: a \cap C_n \neq \emptyset\}: a \in \mathcal A\}\setminus [I]^{<\omega}$ is centered and the finite intersections of this set are cofinite in $I$.
		\begin{theorem}
			Assume CH. Let $\mathbb P$ be a countable forcing poset. Then every infinite countable almost disjoint family may be expanded to $\mathbb P$-indestructible fin-intersecting MAD family $\mathcal A$ such that for $\Vdash_{\mathbb P}``\check{\mathcal A}\text{ is fin-intersecting}''$. In other words, $\mathbb P$-indestructible fin-intersecting MAD families for which forcing with $\mathbb P$ preserves the fin-intersecting property exist generically.
		\end{theorem}
		
		\begin{proof}
			Working in $V$, let $\{a_n: n\in\omega\}$ be an infinite countable almost disjoint family.
			
			Let $((\tau_\gamma, \dot C_\gamma, p_\gamma): \omega\leq \gamma<\omega_1)$ be a listing of all pairs $(\tau, \dot C, p)$ such that:
			
			\begin{itemize}
				\item $\tau$ is a $\mathbb P$-nice name for a subset of $\check \omega$,
				\item $\dot C$ is a $\mathbb P$-nice name for a subset of $\check{(\omega\times [\omega]^{<\omega})}$,
				\item $p \in \mathbb C$,
				\item $p \Vdash \dot C \text{ is a fin sequence on }\omega$.
				
			\end{itemize}
		
			Also, fix $(f_\gamma: \omega\leq\gamma<\omega_1)$ such that $f_\gamma:\omega\rightarrow \gamma$ is bijective for every $\gamma<\omega_1$.

			For $\alpha \in [\omega, \omega_1)$, we recursively define $a_\alpha$ and $\mathcal A_\alpha$ such that, for every $\alpha\in [\omega, \omega_1)$:
			
			\begin{enumerate}[label=\alph*)]
				\item $\mathcal A_\alpha=\{a_\xi: \xi<\alpha\}$ is an almost disjoint family,
				\item for every infinite $\gamma \leq\alpha$ such that $p_\gamma \Vdash`` \dot C_\gamma \text{ is a fin sequence}$'', and for every $F \in [\alpha]^{<\omega}$, if for all $J \in [\alpha]^{<\omega}$: \\
				
				$p_\gamma \Vdash |\{n \in I(\check{\mathcal A|\gamma}, \dot C_\gamma, \check f_\gamma): \forall \xi \in \check F\, (\dot C_\gamma(n)\cap \check a(\xi)\neq \emptyset) \andd  \dot C_\gamma(n)\setminus \bigcup_{\xi \in \check J}{\check a(\xi)}\neq \emptyset\}|=\omega$,\\
				
				then $p_\gamma \Vdash |\{n \in I(\check{\mathcal A|\gamma}, \dot C_\gamma, \check f_\gamma): \forall \xi \in \check F\cup \{\check \alpha\}\,(\dot C_\gamma(n)\cap \check a(\xi)\neq \emptyset) \}|=\omega$.
				
				\item If $p_\alpha \Vdash \tau_\alpha \in \mathcal I^+(\check{\mathcal A|\alpha})$, then $p_\alpha\Vdash |\check a_\alpha\cap \tau_\alpha|=\omega$.
			\end{enumerate}
			
			We show how to construct $a_{\alpha}$ satisfying a), b), c).
			
			First we will define $a^0_\alpha$ satisfying a) and c). Then we will define $a^1_\alpha$ satisfying a) and b). Then we define $a_\alpha=a^0_\alpha\cup a^1_\alpha$. which is then easily seen to satisfy a), b) and c).
			
			Defining $a^0_\alpha$: if the hypothesis of c) fails, just let $a^0_\alpha$ be an infinite set almost disjoint with the infinite countable almost disjoint family $\mathcal A_\alpha$. Now suppose that $p_\alpha \Vdash \tau_\alpha \in \mathcal I^+(\check{\mathcal A|\alpha})$. Let $b_n=a_{f_\alpha(n)}$ for $n \in \omega$. $a^0_\alpha$ will be $\{x_n:n\in\omega\}$, where $x_n$ is defined as follows: we enumerate all pairs $(r, l)$ such that $r\leq p_\gamma$ and $l\in\omega$ as $(r_n, l_n)_{n\in\omega}$. For each $n$, $r_n\Vdash |\tau_n\setminus \bigcup_{i\leq \check n}\check b_i|=\omega$, so there exists $q_n\leq r_n$ and $x_n\in\omega$ such that $x_n\leq l_n$, $x_n\notin \bigcup_{i\leq  n}b_i$ and $q_n\Vdash \check x_n\in\tau_n$. The conclusion is now clear.
			
			Now we define $a^1_\alpha$. Suppose $\{(r, F, \gamma, l): l \in \omega, r\leq p_\gamma, F\in [\alpha]^{<\omega}, \gamma<\alpha, \forall J \in [\alpha]^{<\omega}\,(p_\gamma \Vdash |\{n \in I(\check{\mathcal A|\gamma}, \dot C_\gamma, \check f_\gamma): \forall \xi \in \check F\, (\dot C_\gamma(n)\cap a_\xi\neq \emptyset) \andd \dot C_\gamma(n)\setminus \bigcup_{\xi \in J}{a_\xi}\neq \emptyset\}|=\omega)\}$ is nonempty and enumerate it as $\{(r_m, F_m, \gamma_m, l_m): m \in \omega\}$.
			
			For every $m \in \omega$, there exists $s_m\leq r_m$, $n_m, k_m > l_m$ such that $s_m\Vdash \check n_m \in I(\check{\mathcal A|{\gamma_m}}, \dot C_{\gamma_m}, \check f_{\gamma_m}), \forall \xi \in \check F_m \,\dot C_{\gamma_m}(\check n_m)\cap \check a_\xi\neq \emptyset \text{ and } \check k_m\in \dot C_{\gamma_{m}}(n_m)\setminus \bigcup_{i\leq m} \check a_{\gamma_i}.$
			
			$k_m$ may be picked greater than $l_m$ since $r_m\leq p_{\gamma_m}\Vdash (\text{the }\dot C_{\gamma_m}(n)$'s are pairwise disjoint). Let $a_\alpha=\{k_m: m \in \omega\}$. If the preceding set is empty, just let $a_\alpha$ be an infinite subset of $\omega$ almost disjoint from every $a_\xi$ $(\xi <\alpha)$. This makes $a_{\alpha}$ satisfy b).

			Now let $\mathcal A=\{a_\xi:\xi<\omega_1\}$. By a), $\mathcal A$ is an almost disjoint family. By c), $\mathbb P\Vdash \mathcal A\text{ is MAD}$. We claim that $\mathbb P\Vdash \mathcal A\text{ is fin-intersecting}$ as well. To see that, let $G$ be $\mathbb P$-generic over $V$ and let $C$ be a fin sequence in $V[G]$. If there exists $I$ in $V[G]$ such that $I\in [\omega]^\omega$ and $\bigcup_{n\in I}C_i\in \mathcal I(\mathcal A)$ we are done by Lemma \ref{lemma:idealImpliesFinIntersecting}. Let $G$ be $\mathbb P$-generic over $V$. Suppose that $\mathcal A$ is not MAD in $V[G]$. There exists $\tau \in V^{\mathbb P}$ such that $\tau_G\subseteq \omega$ is infinite and $\tau_G\cap a_\xi$ is finite for every $\xi<\omega_1$. In particular, $\tau_G \in \mathcal I^+(\mathcal A)$. There exists $p \in G$ such that $p\Vdash \tau\subseteq \omega$, $p\Vdash \tau\perp \check{\mathcal A}$ and $p\Vdash \tau \in I^+(\check{\mathcal A})$. There exists $\alpha \in [\omega, \omega_1)$ such that $p_\alpha=p$ and $p_\alpha \Vdash \tau_\alpha=\tau$. In particular, $p_\alpha\Vdash \tau_\alpha \in I^+(\check{\mathcal A|\alpha})$. But then $p_\alpha\Vdash |\tau_\alpha\cap \check a_\alpha|=\omega$, a contradiction.
			
			$\mathcal A$ is indestructibly fin-intersecting: let $G$ be $\mathbb P$-generic over $V$. Fix a fin sequence $C \in V[G]$ such that for all infinite $I\subseteq \omega$ $\bigcup_{n \in I}C_n \in \mathcal I^+(\mathcal A)$. We will show that there exists $I\subseteq \omega$ such that $\{\{n \in I: a\cap C_n\neq \emptyset\}$ is centered. Let $\dot C$ be a name for $C$. Let $p\in G$ be such that $p\Vdash ``\dot C$ is a fin sequence and $\forall I \in [\omega]^\omega \bigcup_{n \in I} \dot C(n) \in \mathcal I^+(\mathcal A)$''. Let $\gamma \in [\omega, \omega_1)$ be such that $p=p_\gamma$ and $p_\gamma\Vdash \dot C=\dot C_\gamma$. Now we will prove the following:
			
			$$p_\gamma\Vdash \forall F \in [\omega_1]^{<\omega}\, \left(\forall \xi \in F\, |\{n \in I(\check{\mathcal A}_\gamma, \dot C_\gamma, \check f_\gamma): a_\xi\cap \dot C_\gamma(n)\neq \emptyset\}|=\omega\right.$$$$\left.\rightarrow \left|\bigcap_{\xi \in F}\{n \in I(\check{\mathcal A}_\gamma, \dot C_\gamma, \check f_\gamma): a_\xi\cap \dot C_\gamma(n)\neq \emptyset\}\right|=\omega\right)$$
			
			(We interpret an intersection with empty domain as $I(\check{\mathcal A}_\gamma, \dot C_\gamma, \check f_\gamma)$). By the preservation of $\omega_1$ this is equivalent to show that for all $F \in [\omega_1]^{<\omega}$:
			
			$$p_\gamma\Vdash \forall \xi \in \check F \,|\{n \in I(\check{\mathcal A}_\gamma, \dot C_\gamma, \check f_\gamma): \check a(\xi)\cap \dot C_\gamma(n)\neq \emptyset\}|=\omega$$$$\rightarrow \left|\bigcap_{\xi \in \check F}\{n \in I(\check{\mathcal A}_\gamma, \dot C_\gamma, \check f_\gamma): \check a(\xi)\cap \dot C_\gamma(n)\neq \emptyset\}\right|=\omega.$$
			
			We will prove the following stronger claim:
			
			$$p_\gamma\Vdash \forall \xi \in \check F\cap \check\gamma \,|\{n \in I(\check{\mathcal A}_\gamma, \dot C_\gamma, \check f_\gamma): \check a(\xi)\cap \dot C_\gamma(n)\neq \emptyset\}|=\omega$$$$\rightarrow \left|\bigcap_{\xi \in \check F}\{n \in I(\check{\mathcal A}_\gamma, \dot C_\gamma, \check f_\gamma): \check a(\xi)\cap \dot C_\gamma(n)\neq \emptyset\}\right|=\omega.$$
			
			Fix $F$. By the definition of $I(\check{\mathcal A}_\gamma, \dot C_\gamma, \check f_\gamma)$ we know that:
			
			$$p_\gamma\Vdash \forall \xi \in \check F\cap \check\gamma \,|\{n \in I(\check{\mathcal A}_\gamma, \dot C_\gamma, \check f_\gamma): \check a(\xi)\cap \dot C_\gamma(n)\neq \emptyset\}|=\omega$$$$\rightarrow \left|I(\check{\mathcal A}_\gamma, \dot C_\gamma, \check f_\gamma)\setminus \bigcap_{\xi \in \check F\cap \check \gamma}\{n \in I(\check{\mathcal A}_\gamma, \dot C_\gamma, \check f_\gamma): \check a(\xi)i\cap \dot C_\gamma(n)\neq \emptyset\}\right|<\omega.$$
			
			In particular, we are done if $F\subseteq \gamma$. If not, let $F\setminus\gamma=\{\alpha_0, \dots, \alpha_{k-1}\}$ such that $\alpha_0<\dots<\alpha_{k-1}$ and let, for $i\leq k$, $F_i=\{\alpha_0\dots, \alpha_{i-1}\}$ (so in particular $F_0=\emptyset$). It suffices to see that for all $i\leq k$:
			
			$$p_\gamma\Vdash \left|\bigcap_{\xi \in \check F_i}\{n \in I(\check{\mathcal A}_\gamma, \dot C_\gamma, \check f_\gamma): \check a(\xi)\cap \dot C_\gamma(n)\neq \emptyset\}\right|=\omega.$$

			Of course, this is true for $F_0=\emptyset$. Suppose $i<k$ and that the above is true for $i$. We show that it is also true for $i+1$. We know that: $p_\gamma\Vdash \left|\{n \in I(\check{\mathcal A}_\gamma, \dot C_\gamma, \check f_\gamma): \forall \xi \in F_i \,\check a(\xi)\cap \dot C_\gamma(n)\neq \emptyset\}\right|=\omega$. Let $\dot L$ be a name such that $p_\gamma\Vdash \dot L=\{n \in I(\check{\mathcal A}_\gamma, \dot C_\gamma, \check f_\gamma): \forall \xi \in F_i \,\check a(\xi)\cap \dot C_\gamma(n)\neq \emptyset\}$. By hypothesis, $p_\gamma\Vdash \bigcup_{n \in \dot L} \dot C(n) \in \mathcal I^+(\mathcal A)$. This implies that if $J \in [\alpha_i]^{<\omega}$, $p_\gamma\Vdash \left|\bigcup_{n \in \dot L} \dot C(n) \setminus \bigcup_{\xi \in \check J} \check a(\xi)\right|=\omega$, which in turn implies that $p_\gamma\Vdash |\{n \in \dot L: \dot C(n) \setminus \bigcup_{\xi \in \check J} \check a(\xi)\neq \emptyset\}|=\omega$. By the definition of $\dot L$ this shows that $p_\gamma\Vdash |\{n \in I(\check{\mathcal A}_\gamma, \dot C_\gamma, \check f_\gamma): \forall \xi \in \check F_i\, \dot C(n)\cap \check a(\xi)\neq \emptyset \andd  \dot C(n) \setminus \bigcup_{\xi \in \check J} \check a(\xi)\neq \emptyset\}|=\omega$. Since $J$ is arbitrary, it follows from (ii) that $p_\gamma\Vdash |\{n \in I(\check{\mathcal A}_\gamma, \dot C_\gamma, \check f_\gamma): \forall \xi \in \check F_{i+1}\, \dot C(n)\cap \check a(\xi)\neq \emptyset|=\omega$. This finishes the proof.
			
		\end{proof}

		\begin{proposition}Assume $\mathcal A$ is an almost disjoint family such that remains MAD and fin-intersecting after adding a Cohen real. Then it also remains MAD and fin-intersecting after adding $\kappa$ Cohen reals with finite support for every $\kappa$.
		\end{proposition}
		
		\begin{proof}Let $\mathbb P$ be a forcing poset that adds $\kappa$ reals and let $G$ be $\mathbb P$-generic over $V$. Let $C$ be a fin sequence and $X$ be an infinite subset of $\omega$, both in $V[G]$. Let $\mathbb Q$ be a (countable) forcing poset that adds a Cohen real. There exists a $\mathbb Q$-generic filter $H$ over $V$ such that $V[H]\subseteq V[G]$ and $C, X\in V[H]$. Then by the hypothesis there exists $a\in\mathcal A$ such that $a\cap X$ is infinite and an infinite $I\in [\omega]^\omega$ in $V[H]$ such that $\{\{n\in I:C_n)\cap a\neq\emptyset\}:a\in\mathcal A\}\setminus[I]^{<\omega}$ is centered. Of course, both things also happen in $V[G]$.
		\end{proof}
	
		Notice that since $\mathfrak s=\omega_1$ after adding any number of Cohen reals the preceding result is non-trivial (see Theorem \ref{<simpliesfinintersecting}).
		\section{A fin-indestructible MAD family in the Random Model}
		
		In this section we show that assuming CH, every countable almost disjoint family can be extended to a fin-intersecting MAD family which remains a fin-intersecting MAD family after adding an arbitrary quantity of Random reals (with finite supports). This construction has many similarities to the previous one, but uses some different techniques to deal with the fact that the forcing notions we are dealing with now are not countable. Some of the techniques used are similar to the construction of indestructible MAD families for some proper forcing notion (see \cite[Lemma III.1.]{hruvsak2001life}).
		\begin{theorem}
			Assume CH. Let $\mathbb P$ be a proper forcing poset of size $\leq \omega_1$ such that for every $\dot f \in V^{\mathbb P}$ and $p \in \mathbb P$, if $p\Vdash \dot f \in \omega^\omega$ then there exists $h \in \omega^\omega$ such that $p\Vdash \dot f \leq^* h$.  Then every infinite countable almost disjoint family may be expanded to a $\mathbb P$-indestructible fin-intersecting MAD family $\mathcal A$ such that for $\Vdash_{\mathbb P}``\check{\mathcal A}\text{ is fin-intersecting}''$. In other words, $\mathbb P$-indestructible fin-intersecting MAD families for which forcing with $\mathbb P$ preserves the fin-intersecting property exist generically.
		\end{theorem}
		\begin{proof}
			Let $\mathcal A'$ be a given infinite countable almost disjoint family and write $\mathcal A'=\{a_n: n \in \omega\}$ injectively.  By properness, there exists a family of triples $((p_\alpha, \tau_\alpha, \dot C_\alpha): \omega\leq \alpha<\omega_1)$ such that whenever $p\in \mathbb P$, $\tau \in V^{\mathbb P}$, $\dot C \in V^{\mathbb P}$ and $p\Vdash \tau\subseteq \omega$ and $p\Vdash \dot C \text{ is a fin sequence}$, there exists $\alpha \in [\omega, \mathfrak c)$ such that $p_\alpha\leq p$, $p_\alpha \Vdash \dot C=\dot C_\alpha$ and $p_\alpha \Vdash \tau=\tau_\alpha$.
			
			Fix a family $(f_\alpha: \alpha \in [\omega, \omega_1))$ such that for each $\alpha \in [\omega, \omega_1)$, $f_\alpha:\alpha\rightarrow \omega$ is bijective.
			
			Enumerate $\mathcal A'=\{a_n: n \in \omega\}$. For $\alpha \in [\omega, \omega_1)$, we recursively define $a_\alpha$ and $\mathcal A_\alpha$ such that, for every $\alpha\in [\omega, \omega_1)$:
			
			\begin{enumerate}[label=\alph*)]
				\item $\mathcal A_\alpha=\{a_\xi: \xi<\alpha\}$ is an almost disjoint family,
				\item for every infinite $\gamma \leq\alpha$ such that $p_\gamma \Vdash`` \dot C_\gamma \text{ is a fin sequence}$'', and for every $F \in [\alpha]^{<\omega}$, if for all $J \in [\alpha]^{<\omega}$: \\
				
				$p_\gamma \Vdash |\{n \in I(\check{\mathcal A|\gamma}, \dot C_\gamma, \check f_\gamma): \forall \xi \in \check F\, (\dot C_\gamma(n)\cap \check a(\xi)\neq \emptyset) \andd  \dot C_\gamma(n)\setminus \bigcup_{\xi \in \check J}{\check a(\xi)}\neq \emptyset\}|=\omega$,\\
				
				then $p_\gamma \Vdash |\{n \in I(\check{\mathcal A|\gamma}, \dot C_\gamma, \check f_\gamma): \forall \xi \in \check F\cup \{\check \alpha\}\,(\dot C_\gamma(n)\cap \check a(\xi)\neq \emptyset) \}|=\omega$.
				
				\item If $p_\alpha \Vdash \tau_\alpha \in \mathcal I^+(\check{\mathcal A|\alpha})$, then $p_\alpha\Vdash |\check a_\alpha\cap \tau_\alpha|=\omega$.
			\end{enumerate}
			
			We show how to construct $a_{\alpha}$ satisfying a), b), c).
			
			First we will define $a^0_\alpha$ satisfying a) and c). Then we will define $a^1_\alpha$ satisfying a) and b). Then we define $a_\alpha=a^0_\alpha\cup a^1_\alpha$ and this will satisfy a), b) and c).
			
			Defining $a^0_\alpha$: If the hypothesis of c) fails, just let $a^0_\alpha$ be an infinite set almost disjoint to the infinite countable almost disjoint family $\{a_{\xi}:\xi<\alpha\}$. Now suppose that $p_\alpha \Vdash \tau_\alpha \in \mathcal I^+(\check{\mathcal A|\alpha})$. Let $b_n=a_{f_\alpha(n)}$ for $n \in \omega$. Let $\rho$ be a name such that $p_\alpha \Vdash \rho:\omega\rightarrow \omega \text{ is strictly increasing}$ and $p_\alpha \Vdash \forall n \in \omega \rho(n)\in \tau_\alpha\setminus \bigcup_{i<n} \check b(i)$.  Let $h:\omega\rightarrow \omega$ be such that $p_\alpha \Vdash \rho<^* \check h$. Let $a^0_\alpha=\bigcup_{n \in \omega} h(n)\setminus \bigcup_{i<n}b_\alpha(i)$:
			
			$p_\alpha \Vdash |\check a^0_\alpha\cap \tau_\alpha|=\omega$ (which implies that $a^0_\alpha$ is infinite): Suppose by contradiction that this is false. Then there exists exists $n \in \omega$ and $q\leq p_\alpha$ such that $q\Vdash \check a^0_\alpha\cap \tau_\alpha \subseteq \rho(\check n)\andd \rho(\check n)<\check h(\check n)$, but then $q\Vdash \rho(\check n)\in \check h(\check n)\setminus \bigcup_{i<\check n}\check b_\alpha(i)\subseteq \check a_{\alpha}^0$ and $q\Vdash \rho(\check n) \in \tau_\alpha$, a contradiction. In particular, this shows that $a^0_\alpha$ is infinite.
			
			$a^0_\alpha$ is clearly almost disjoint from each $a_{\xi}$ for $\xi<\alpha$.

			Now we define $a^1_{\alpha}$. 
			If the hypothesis of b) fails, just let $a_\alpha^1\in [\omega]^\omega$ be almost disjoint with $a_\alpha'$ for every $\xi<\alpha$. If it does not fail, we proceed as follows:

			For each infinite $\gamma\leq \alpha$ such that $p_\gamma\Vdash \dot C_\gamma \text{ is a fin sequence}$, let $\dot K_\gamma$ denote a name such that $p_\gamma\Vdash \dot K_{\gamma}: [\alpha]^{<\omega}\times [\alpha]^{<\omega}\rightarrow \mathcal P(\omega)$ and: $$p_\gamma\Vdash \forall F, J \in [\check \alpha]^{<\omega}\, \dot K_{\gamma}(F, J)$$$$=\left\{n \in I(\check{\mathcal A|\gamma}, \dot C_\gamma, \check f_\gamma): \forall \xi \in \check F\, (\dot C_\gamma(n)\cap \check a(\xi)\neq \emptyset) \andd  \dot C_\gamma(n)\setminus \bigcup_{\xi \in J}{\check a(\xi)}\neq \emptyset\right\}.$$
			
			Now let:
			
			$$\mathcal F_{\gamma}=\left\{F \in [\alpha]^{<\omega}: \forall J \in [\alpha]^{<\omega}\, p_\gamma \Vdash |\dot K_{\gamma, \alpha}(\check F, \check J)|=\omega\right\}$$
			
			
			Fix $\gamma$. Suppose $\mathcal F_{\gamma}\neq \emptyset$. Enumerate $\mathcal F_{\gamma}\times \omega=\{(F_n^{\gamma}, l_n^{\gamma}): n \in \omega\}$. Now let $\rho_{\gamma}$ be a name such that $p_\gamma \Vdash \rho_{\gamma}:\omega\rightarrow \omega$	and $p_\gamma \Vdash \forall n \in \omega\, \rho_{\gamma}(n)\in \dot K_{\gamma}(F_n^\gamma, \check{f}_\alpha[n])\setminus \check l_n^\gamma$. Let $\bar \rho_{\gamma}$ be a name such that $p_\gamma \Vdash \bar \rho_\gamma:\omega\rightarrow \omega$ and $p_\gamma\Vdash \forall n \in \omega \, \bar \rho_{\gamma}(n)=\max \left(\dot C_\gamma(\rho_{\gamma}(n))\setminus \bigcup_{i<n}\check a_{\check f_\alpha(i)}\right)$. Fix  $h_\gamma:\omega\rightarrow \omega$ such that $p_\gamma\Vdash \bar \rho_\gamma<^*\check h_\gamma$.
			
			Let $U=\{\gamma\leq \alpha: \gamma\geq \omega\andd p_\gamma\Vdash ``\dot C_\gamma \text{ is a fin sequence''}\andd\mathcal F_{\gamma, \alpha}\neq \emptyset\}$. If $U$ is empty, there is nothing to do, so just let $a^1_\alpha$ be infinite and almost disjoint with every $a_\gamma$ for $\gamma<\alpha$. It $U$ is not empty let $h:\omega\rightarrow \omega$ be such that $h\geq^* h_{\gamma}$ for every $\gamma\in U$. Let $a_\alpha^1=\bigcup_{n \in \omega} h(n)\setminus \bigcup_{n \in \omega} a_{f_\alpha(n)}$. Of course $a_\alpha^1$ is almost disjoint with $a_\gamma$ for every $\gamma<\alpha$. Now we verify b) (which implies, in particular, that $a_\alpha$ is infinite).
			
			Suppose $\gamma\leq \alpha$ and $F$ satisfy the hypothesis of b). So $F \in \mathcal F_{\gamma}$. We must see that: $$p_\gamma \Vdash |\{n \in I(\check{\mathcal A|\gamma}, \dot C_\gamma, \check f_\gamma): \forall \xi \in \check F\cup \{\check \alpha\}\,(\dot C_\gamma(n)\cap \check a(\xi)\neq \emptyset) \}|=\omega.$$
			
			Suppose this is false. So there exists $r\leq p_\gamma$ and a natural number $l$ such that:
			
			$$r \Vdash \{n \in I(\check{\mathcal A|\gamma}, \dot C_\gamma, \check f_\gamma): \forall \xi \in \check F\cup \{\check \alpha\}\,(\dot C_\gamma(n)\cap \check a(\xi)\neq \emptyset) \}\subseteq \check l $$$$\andd \forall n\geq \check l\, \bar \rho_\gamma(n)<\check h_\gamma(n).$$
			
			There exists $N\geq l$ such that $h(m)>h_{\gamma}(m)$ for every $m\geq N$. Let $n\geq N$ be such that $l_n^{\gamma}\geq l$ and $F_n^{\gamma}=F$. Then $r\Vdash \rho_{\gamma}(\check n) \in \dot K_{\gamma}(F_n, \check f_\alpha[\check n])\setminus \check l_n^\gamma$, so $r\Vdash\rho_\gamma(\check n)\geq\check l$. In order to get a contradiction, it suffices to see that $r\Vdash \dot C_\gamma(\rho_{\gamma}(\check n))\cap \check a_\alpha\neq \emptyset$. But this is true since $r\Vdash \bar \rho_{\gamma} (\check n) \in \dot C_\gamma(\rho_{\gamma}(\check n))$ and $r\Vdash \bar \rho_{\gamma} (\check n) \in \check h_{\gamma}(\check n)\setminus \bigcup_{i<\check n}\check a\circ \check f_\beta(i)\subseteq \check h_{\alpha}(\check n)\setminus \bigcup_{i<\check n}\check a\circ \check f_\alpha(i)\subseteq \check a_{\alpha}^1$.\\
			
			This finishes the construction. The verification that $\mathcal A=\{a_\xi: \xi<\omega_1\}$ has the required properties is very similar to the one in the construction in the previous section and left to the reader

		\end{proof}

		\begin{corollary}Assume CH and let $\mathbb P_{\kappa}$ be the standard forcing notion for adding $\kappa$ random reals. Then every infinite countable almost disjoint family may be expanded to $\mathbb P$-indestructible MAD family $\mathcal A$ such that for $\Vdash_{\mathbb P_\kappa}``\check{\mathcal A}\text{ is fin-intersecting}''$. In particular, it remains pseudocompact.
		\end{corollary}
		
		\begin{proof}Let $\mathcal A'$ be an infinite countable almost disjoint family. Let $\mathbb P$ be the forcing notion for adding $\omega$ random reals. By the previous theorem there exists a MAD family $\mathcal A$ containing $\mathcal A'$ which remains MAD and is fin-intersecting after forcing with $\mathbb P$. We argue that the same happens with $\mathbb P_\kappa$.
			
			It remains MAD: let $G$ be a $\mathbb P_{\kappa}$-generic filter over $V$ and $x \in [\omega]^\omega\cap V[G]$. It is well known that there exists a $\mathbb P$-generic filter $H$ such that $x \in V[H]\subseteq V[G]$, thus, there exists $a \in \mathcal A$ such that $a\cap x$ is infinite.
			
			It is fin-intersecting: let $G$ be a $\mathbb P_{\kappa}$-generic filter over $V$ and $C\in  V[G]$ be a fin sequence. Since $C$ can be seen as a real, there exists a $\mathbb P$-generic filter $H$ such that $C \in V[H]\subseteq V[G]$. Assume that, $V[G]\models \forall I \in [\omega]^\omega \bigcup_{n \in I}C_n \in \mathcal I^+(\mathcal A)$. By downards absoluteness, $V[H]\models \forall I \in [\omega]^\omega \bigcup_{n \in I}C_n \in \mathcal I^+(\mathcal A)$, so $V[H]\models \exists I \in [\omega]^\omega \{\{n \in I: a\cap C_n\neq \emptyset\}: a \in \mathcal A\}\setminus [I]^{<\omega} \text{ is centered}$. Now, by upwards absoluteness, $V[G]\models \exists I \in [\omega]^\omega \{\{n \in I: a\cap C_n\neq \emptyset\}: a \in \mathcal A\}\setminus [I]^{<\omega} \text{ is centered}$. This completes the proof.
		\end{proof}

		It is easy to merge the Cohen and Random constructions into a single construction: just enumerate all the names as in the two constructions and instead of uniting just two sets to obtain $a_\alpha$, take the union of four sets, one of them satisfying b) in the Cohen Construction, other c) in the Cohen construction, the third satisfying b) in the Random Construction, and the fourth d) in the Cohen construction. With this modification we obtain the following:
		
		\begin{corollary}Assume CH. Then fin-indestructible MAD families which remain so by adding arbitrarily many Cohen reals or Random reals with finite supports exist generically.
		\end{corollary}
	
		\section{Acknowledgements}
	   The authors thank professors Michael Hru\v{s}\'{a}k and Paul Szeptycki for suggestions. The second author also thanks Professor Artur Tomita for early discussions on this material.
		\bibliographystyle{plain}
		\bibliography{bibliografia}{}
	\end{document}